# Influence of a function's coefficients and feedback of the mathematical work when reading a graph in an online assessment system

Jorge Gaona[1], Rodrigo Hernández[2], Felipe Guevara[3] y Víctor Bravo[2].

**Abstract**

*This paper shows the results of an experiment applied to 170 students from two Chilean universities who solve a task about reading a graph of an affine function in an online assessment environment where the parameters (coefficients of the graphed affine function) are randomly defined from an ad-hoc algorithm, with automatic correction and automatic feedback. We distinguish two versions: one of them with integer coefficients and the other one with decimal coefficients in the affine function. We observed that the nature of the coefficients impacts the mathematical work used by the students, where we again focus on two of them: by direct estimation from the graph or by calculating the equation of the line. On the other hand, a feedback oriented towards the "estimation" strategy influences the mathematical work used by the students, even though a non-negligible group persists in the "calculating" strategy, which is partly explained by the perception of each of the strategies.*

*Digital technology, task design, mathematical activity, distance education and online learning, post-secondary education*

## 1 Introduction

To analyze data or information, to estimate quantifiable amounts of products, events, or information, and to process and interpret the sense of the information are part of the skills required in the XXI century in jobs related to STEAM disciplines (Jang, 2016, p.293). One of the ways to represent data is to work with graphs, hence, its construction, reading, comprehension and interpretation are abilities that students must acquire during their education (Gravemeijer et al., 2017; Hoyles et al., 2010).

The research referring to the use of this type of information representation is abundant, some of it treating the comprehension of a function's graph (Mudaly & Rampersad, 2010), the derivative (Borji et al., 2018; Ergül, 2018; Graham & Sharp, 1999; Haciomeroglu et al., 2010; Mudaly & Rampersad, 2010) or the area under the curve in an integral. (Klein et al., 2019). Other works focus on mathematical and physical contexts (Eshach, 2014; Planinic et al., 2012; Woolnough, 2000) or other extra-mathematical contexts. (Planinic et al., 2013; Roth & Bowen, 2001; Roth & Hwang, 2006). The interpretation of graphs in the resolution of problems has also been part of the treated topics (Jackson et al., 1993), and so has been the relationship between the language and the interpretation of the graphs (Radford, 2009; Smit et al., 2016).

---
[1] Departamento de Ciencias, Universidad Metropolitana de Ciencias de la Educación, Chile.
[2] Facultad de Ingeniería y Ciencias, Universidad Adolfo Ibáñez, Chile.
[3] Facultad de Ingeniería, Universidad de Atacama, Chile

Other research works study the assessment of tasks involving graphs, for example, comparing the construction of a graph versus the election of one in a multiple-choice question, where differences can be seen in the results obtained by the students (Berg & Smith, 1994; Craig & Philip, 1994), how do strategies change when modifying some didactic variables, e.g., scientific data o mathematical words, or when it is asked for a value outside the graph's range (Bragdon et al., 2019).

To interpret and fully understand a graph it is essential to properly read it. The research about the mentioned process is rather scarce. There is, for example, the research of Ludewig et al. (2020) who defines the reading of a graph as: *The capacity to extract and fluidly use its information.* In his research work he established that the estimation, in a numeric line, is one of the four significant predictors concerning the student's performance.

If we consider the previously mentioned definition of graph reading, it is our concern to know if there are different ways to extract the information provided by the graph when working with affine functions with some variations on the coefficients and when, also, there is an online assessment system that mediates the work, validates the student´s answer and gives feedback. More precisely, we intend to inquire what is being measured when asking the students to estimate the image of a given value. More precisely, the questions that lead this investigation are as follows:

- What is the mathematical work that students do when reading a graph if the coefficients that describe an affine function change between whole numbers and decimals?
- Does a directly value-reading oriented feedback have any influence on the mathematical work made by the students?
- Which is the effectiveness of the work made by the students? (effectiveness measured as the number of correct responses over the number of total responses)

The current research was applied to 2 universities at Chile, one of them is a public university and the other one is private. The questionnaire was applied on the second semester of year 2020 to students who possessed previous knowledge of functions. A total of 170 students participated, 80 from the private university, and 90 was provide from the public one. This was made under the context of on-line classes due to the SarsCov-2 (Covid-19) pandemic.

## 2 Theoretical framework

The theoretical approximation of this work is based in three main ideas. It is defined, firstly, the assessment, particularly the learning assessment and the designed tools for it. Technology is then defined and so is the role it plays in this process, its potential and limitations. Finally, the mathematical working space is described, specifically when estimating the image of an element in the graph of the corresponding function.

### 2.1 Assessment

Assessment is a wide concept, and it accepts many classifications, in which are included the purpose with which it is done, the function it fulfills, the evaluators or the moment in which it is done, among many others. To our ends, we distinguish the role it plays, establishing a difference between learning assessment and assessment for learning. The firstly mentioned is the process that seeks to establish, through a precise tool, whether a skill or knowledge has been achieved or not. The second one is the process of searching and interpreting proofs so that the students and their teachers can use them to decide in which learning stage they are, where should they go, and which is the best way to get there (Wiliam, 2011). Regarding the tools of assessment, it is considered that they can be useful for both learning assessment and assessment for learning (Bennett, 2011), hence, its discussion will consider its purpose.

The previously described categories understand, implicitly, the idea of evaluation as a kind of dialogue between the evaluator and the evaluated. This dialogue implies the delivery of information in both directions, being the one given to the evaluated of particular importance. This information is called feedback (Hattie & Timperley, 2007). It is possible to distinguish several types of feedback in a given task: At a task level, at a process level and self-regulation. (Hattie, 2017)

## 2.2 Tasks mediated by digital artefacts

We will also consider tasks mediated by digital artefacts, understanding an artefact as a human construction with a specific purpose and that a subject can transform into an instrument via usage schemes (Rabardel, 1995). To take in consideration the specificity of software used for doing, teaching, and learning math, digital artefacts will be considered as bearers of what Radford names historic intelligence (Radford, 2014), which consists in the epistemology with which software try to be coherent, or to have as a reference when a subject interacts with them. In contrast to the previous concept, the idea of epistemological relativity appears, which is the intrinsic inability that computer programs have for being entirely loyal to those concepts that they attempt to represent (Balacheff, 2000).

On the other hand, the concepts of 'Pragmatic value' and 'Epistemic value' will be taken in consideration when discussing technology-mediated tasks (Artigue, 2002). Pragmatic value is understood as that which technology makes possible or allows to do in a more efficient and effective way. The epistemic value, however, refers to a technology-mediated task's potential for contributing to the comprehension of mathematical concepts. These elements form a virtual space on which a task can be designed (Artigue, 2002). Pragmatic value is understood as that which technology makes possible or allows to do in a more efficient and effective way. The epistemic value, however, refers to a technology-mediated task's potential for contributing to the comprehension of mathematical concepts. These elements form a virtual space on which a task can be designed.

### 2.2.1 The Task

From Gaona's work (2018, p. 85) a task is decomposed in three components: the type of task, the mathematical objects involved and the context. To define a specific task, there can be a series of variables in each of these components, which is shown in Figure 1.

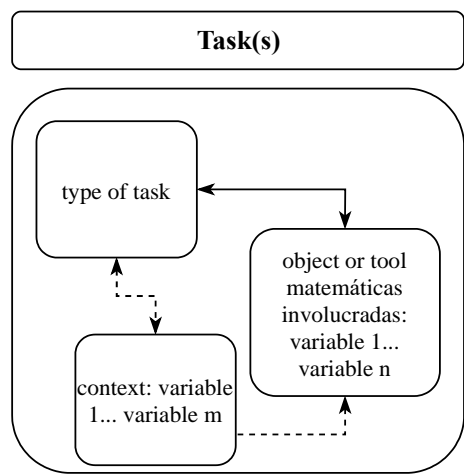

Figure 1: components of a task (own elaboration)

On the other hand, there are *mathematical tools or objects* associated to the *type of task.* This object is observable through its semiotic representation (Duval, 1993). Accordingly, for example, if the type of task is *to estimate the image of a value from a function* (Figure 4), some of the variables are: the type of function, the used representations (graphs, algebraical, tabular), among others. Also, there are other variables involved within each representation. If it is a graphic representation, the line graph, the grid, the graduality and the relationship between each (e.g., if the line passes through certain values or quadrants) would be other variables and so on. The answer (if there is one) could also be (or contain) a mathematical object. So, didactic variables can be identified in the answer just as they can be identified in the mathematical object enunciated. Continuing with the example, is the answer a number or a point on the graph? Do the numbers defining these objects contain integer values, decimal values, or some other type of value from another set? Each of these questions are variables of the task. Duval, 1993). Accordingly, for example, if the type of task is *to estimate the image of a value from a function* (Figure 4), some of the variables are: the type of function, the used representations (graphs, algebraical, tabular), among others. Also, there are other variables involved within each representation. If it is a graphic representation, the line graph, the grid, the graduality and the relation between each (e.g., if the line passes through certain values or quadrants) would be other variables and so on. The answer (if there is one) could also be (or contain) a mathematical object. So, didactic variables can be identified in the answer just as they can be identified in the mathematical object enunciated. Continuing with the example, is the answer a number or a point on the graph? Do the numbers defining these objects contain integer values, decimal values, or some other type of value from another set? Each of these questions are variables of the task.

It seems that the *context* is linked by a dashed line because there can be tasks that require no context and those that do. If there is one, the mathematical object or the type of task at stake may be influenced by the context. For example, if the task is still *to estimate the image of a value from a function* and the context is the cost as a function of quantity bought, this would imply that the domain is both positive and discrete, the numbers could have higher values, among other possible changes.

### 2.2.2 Digital Support

Finally, we define the support on which the task is worked. The task could be delivered in different formats and go along with diverse kinds of support. In this work we will focus on tasks supported by a specific digital artefact: An online assessment system. In Gaona (2020) the artefact brakes down into four components: statement, input system, validation system and feedback.

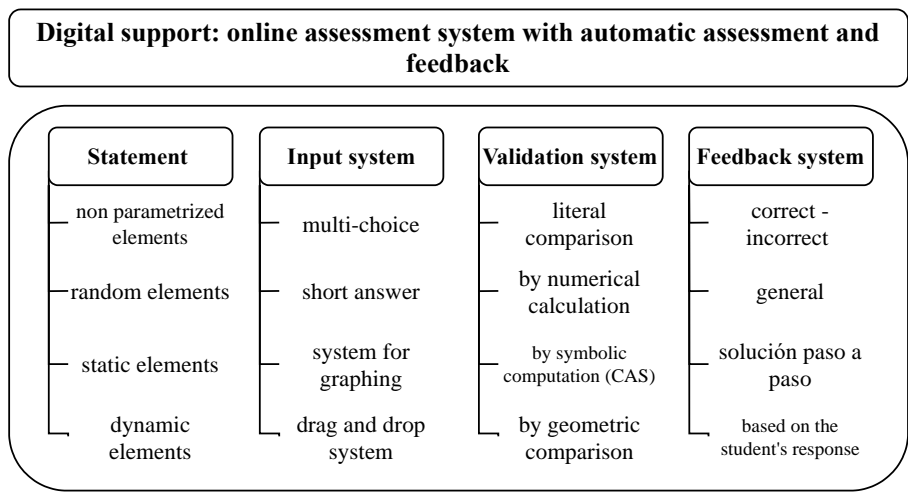

Figure 2: Components of an artefact-task in an online assessment system (own elaboration)

The *statement* shows the proposed task. It can contain fixed or parameterized elements. The last case mentioned implies that each student faces a similarly structured question containing different values. The statement can also contain static (like images) or dynamic elements (like applets with sliders).

The Input System allows students to enter the answer. There are several formats, like using plain text, equation editors, handwrite recognition system, sliders, or multiple-choice answers, among others.

For the validation system we consider the applications that have a geometric processor or a symbolic calculation system. When the answer is not in multiple-choice format it is essential that the system meets certain level of sophistication for the validation. For example, if the answer is an algebraical expression or a number, the system must recognize if two expressions are equivalent. It can also be important that some of the object's characteristics is recognized, like if an expression has been factorized, if a value is within a given range or a matrix fulfills certain conditions regarding its coefficients or regarding the entire matrix.

Finally, the *feedback system* delivers information after the system has already validated the student's answer. Depending on the characteristics of the system, the feedback can indicate if the answer is correct or not, give a step-by-step feedback according to the statement, or a feedback according to the given answer, among others.

The online assessment system's capacity to automatically correct an answer could be considered as a part of its pragmatic value. However, this capacity gives no insight on its epistemic value. This will depend on how the mathematical objects are represented in the statement and what is the task asked to the students, the format in which the answer must be entered, the validation process and what is the feedback given to the student. This can all affect the students, mathematical work.

## 2.3 Mathematical Working Space

A Mathematical Working Space (MWS) is conceived as an abstract structure that allows the study of epistemological and cognitive aspects when individuals solve problems in a specific domain, such as algebra, analysis, geometry, or probabilities, among others (Kuzniak, Tanguay, et al., 2016).

The epistemological and cognitive aspects are articulated through semiotic, instrumental, and discursive genesis (See Figure 3). The word 'genesis' is here used in a wide sense, and it refers, not just to the beginning of a process, but also to its development and interaction between the cognitive and epistemological planes poles.

The *semiotic genesis* links the visualization process at the cognitive plane with the *representamen* at the epistemological one. This genesis could begin due to the sign at the *representamen* which is decodified by the subject through visualization. It could also begin by the subject's codification in which a sign is produced.

The *instrumental genesis* links the construction process at the cognitive plane with the artefacts pole. When working con material, IT, or symbolic tools it is composed of two processes: Instrumentalization and instrumentation (Rabardel, 1995). The first process understands the emergency and evolution of usage schemes of the artefact and the utilization of possibilities that it offers. The second process starts at the subject, and it is relative to the emergency and the evolution of usages and instrumented actions schemes, its constitution, functioning, coordination, combination, inclusion, and assimilation of new artefacts to schemes already constituted. The mathematical work could be considered rutinary if it is not connected with the validation and justification of the artefacts.

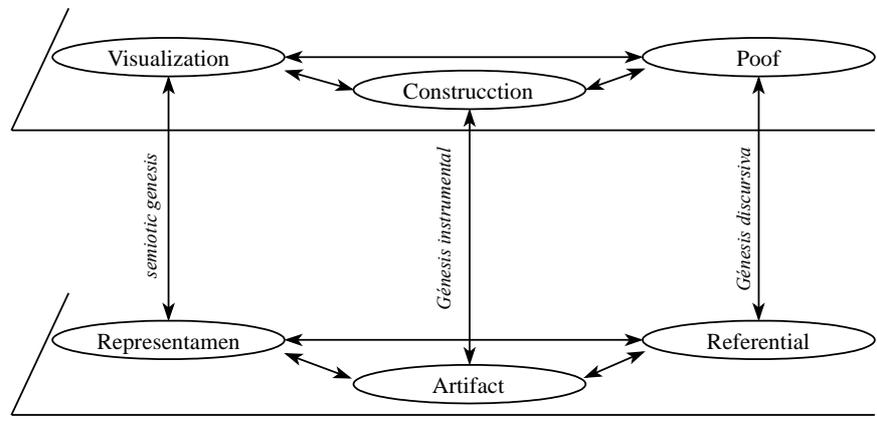

Figure 3: Mathematical Working Space (Kuzniak y Richard, 2014)

At last, discursive genesis connects the proof process with the referential pole on the epistemological plane and it is associated to the deductive reasoning process through theorems and properties. In this last case, the focus is set on the properties and theorems, which is why it is being considered the reasonings that goes further than the visual or instrumental ones, but that can be unchained by them. (Kuzniak, Tanguay y Elia, 2016, p. 727).

When it is not possible to distinguish which genesis is being privileged, which happens frequently, the mathematical work can be characterized through the connection of two geneses, considering some of the three vertical planes (Coutat y Richard, 2011): semiotic-instrumental, semiotic-discursive, or instrumental-discursive.

In this article's introduction, graph reading was defined as the capacity of extracting and fluidly using the graph's information. The way to extract this information may vary depending on the task. In our case, we work with a task that consists in estimating the image of a value existing in the domain. For this specific task we consider two possible strategies: 1) estimating the value from a direct reading and 2) calculating the linear equation from the visualization of two points and then evaluating the equation in the pre-image stated.

In the first case, the estimation is made according to three processes defined by Pizarro, Albarracín and Gorgorió (2018): assigning a value, to perceptively execute the task, and linking the perception with the previous knowledge or with the mental image of auxiliar object. In terms of the MWS, the graph's elements are used as semiotic tools, the individual makes a proportional calculus taking as referent the magnitude of the height of the rectangle from de grid's graduation.

For the second case, an estimation is also made to obtain two ordered pairs to find the linear equation ($(y-y_1) = [(y_2-y_1)/(x_2-x_1)] *(x-x_1)$), but in this case, it is the linear equation, as a symbolic artefact, that leads the work. It guides the convenient points to visualize for using the equation and, consequentially, a work on the semiotic-instrumental plane unfolds.

# 3 Methodology

## 3.1 Materials

To answer the research questions, with the Wiris Software, a task was designed on the platform 'Moodle' which consisted in asking the image of a value from a function represented in shape, as shown in Figure 4.

From Gaona's classification of a task in an online assessment system (2020), said task is broken down into four elements:

1. Statement: The statement is made of the type of task: 'to estimate the image of a value', the value on which the image is asked, and function shown as its graphic shape. In the questions there are (random) parametrized elements defined by an algorithm, so that each student faces a question with the same structure but with different values. The random parameters are:

   o The image asked to evaluate, which is an integer number between –8 and 8 (excluding 0).
   o The coefficients defining the function shown: they can be integer or decimal random numbers, that is, if the linear equation's form is $y=mx+b$ (in both cases), then in the question of Figure 4(a) parameters m and b are integer numbers, and in the Figure 4(b) are decimal numbers with one decimal.
   o The elements of the cartesian plane: center, width, and height. These were defined from an algorithm depending on the parameters of the affine function.

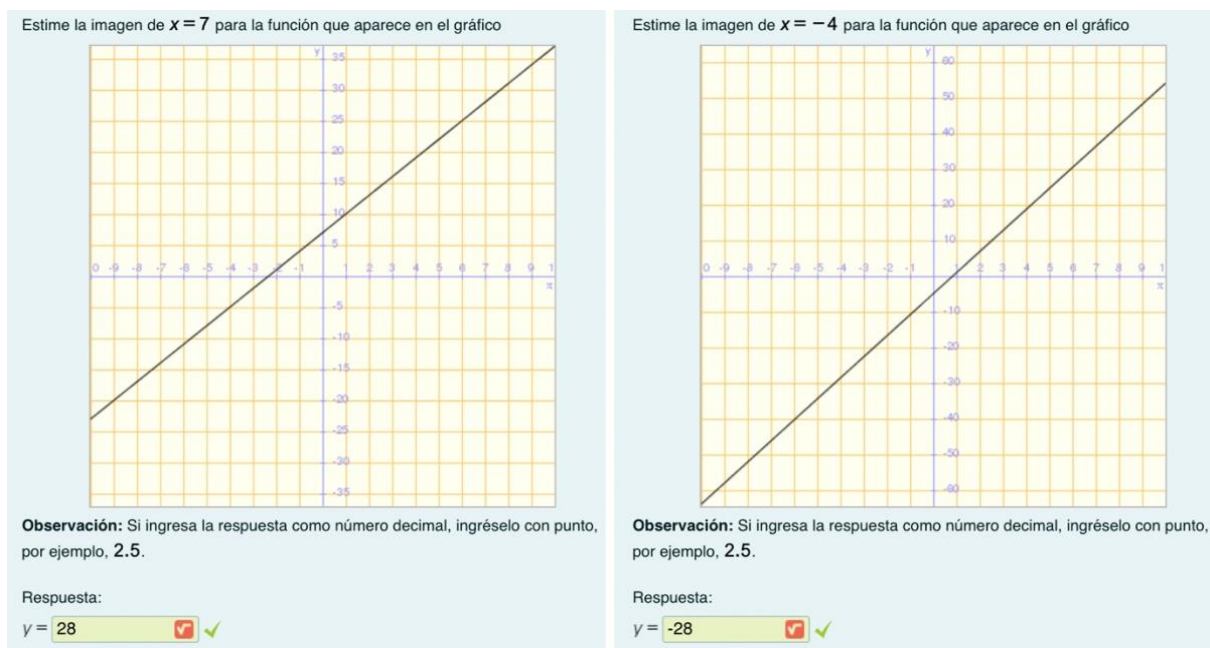

Figure 4: Statements of the task that students had to solve. Left graph (a): the coefficients that define the line are integer numbers, the grid is 1 in 1 on x-axis and 5 in 5 on y-axis. Right graph (b): coefficients defining the line are decimal numbers, the grid is 1 in 1 on x-axis and 10 in 10 on y-axis.

2. Input System: Students have a space to enter a plane text type of answer, but where they can also use an editor to enter fractions or other type of mathematical symbols.
3. Validation System: The platform counts with a symbolic-calculus system that allows to determine if an answer is within the correct range defined on the algorithm through a condition that considers a range within the correct answers (it will be analyzed later).
4. Feedback: Once the student answers the question, the system indicates if the answer is correct or not. Also, it shows a step-by-step solution where the proposed strategy is to draw a vertical line and then a horizontal one for later estimate the value of the image (see Figure 5).

It can be seen that the task generates meaning over the estimation through two mechanisms: the validation of the value in a range, followed by a feedback where this range is explained, and the way for solving the task. For example, in the first case, it is indicated that: "We can observe that the value is between the horizontal line y=-30 and the horizontal line that passes through the center of the stretch between –30 and –20. Then we can assure that the image asked is a value within the range ]-30, -25[". In the second case, it says: *"We can see that its value is between the center and the line y=30. As it is close to the center, we can assure that it is within the third quarter of the stretch between 25 and 30. So, we can assure that the image asked is some value within the interval ]26.25, 28.75["*. Also, a space was provided for the students to write the strategy used to answer.

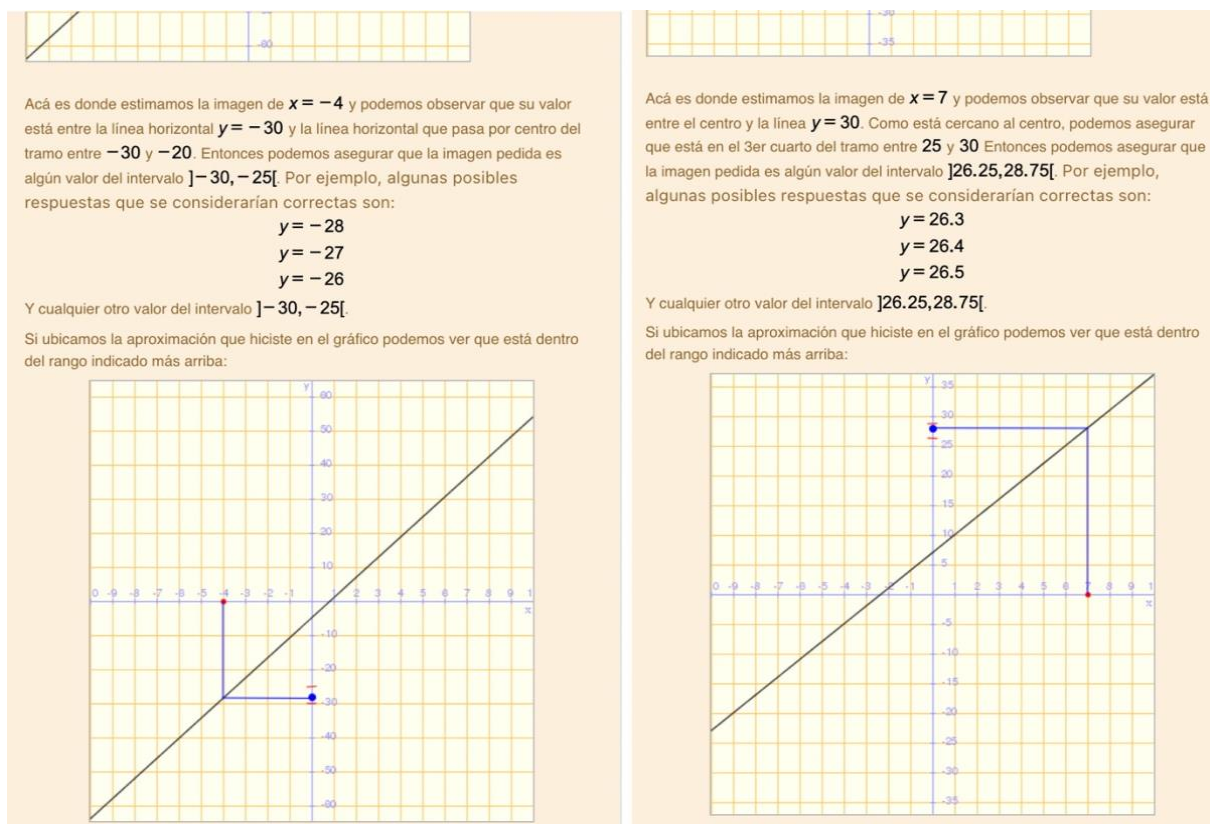

Figure 5: Extract from the feedback of the tasks. On the left is the one in Fig. 4(a) and on the right is the one in Fig. 4(b).

## 3.2 Programming the Ambiguity

An interesting point to highlight is the design of the task and particularly the algorithm that defines it, which can be analyzed from a mathematical, didactic, and instrumental perspectives.

The task proposes a graph of a function of the form *y=mx+b*, where *m* and *b* are random parameters:
- First Case:
  - *m=random(1,5)\*random {-1,1}*
  - *b=random(1,10)\*random {-1,1}*

In this case, lines with integer coefficients between –10 and 10 are produced, without considering 0.

- Second Case:
  - *m=random(1,5)\*random {-1,1} + random(1,9)/10.0*
  - *b=random(1,10)\*random {-1,1} + random(1,9)/10.0*

In this case, lines with decimal coefficients between –10 and 10 are produced, without considering 0.

Then, the cartesian plane where the graph will be drawn is defined, specifically, the center in the point (0,0); the visible distance on the x-axis, between –10 and 10; the height, depending on the higher value achievable when evaluating *y=mx+b* between –10 and 10; the stretch between lines on the x-axis that goes 1 in 1 and the stretch at the y-axis that depends of the previously defined height, and that could go 2 in 2, 5 in 5, 10 in 10, or 20 in 20 (See Figure 4). Finally, the value of the pre-image that will be evaluated is defined, which is a random number between –8 and 8 (without considering 0). From these elements, we define the margin which is considered as correct. According to this, we emphasize in the fact that the measure, in terms of the length of the rectangle, is variable: it can take the value of 1x2, 1x5, 1x10, or 1x20, and the 'how it looks' depends on the height of the graph and of the 400 pixels over which the cartesian plane is drawn on screen.

The task, being about visualization, must consider these elements to measure the rectangle's size in pixels that will see who answers the task. Based on this, the algorithm must measure the height in which the image is found and give feedback to the students estimating how 'close' or how 'far' is the image from the horizontal lines. For example, in Figure 4(a) it must be measured how close is the image of x=-4 to the line y=-30. Generally, it was defined how big or small was the height of the rectangle over which the student had to estimate the value, the distance from the image to the line by separating the number of ranges according to the rectangle's height. From these parameters, it was defined, in the algorithm, the range considered as 'correct´.

Some of the decision that must be made in this process are rather ambiguous, because, as a last resort, we must do it depending on whether the rectangle observed is taller or shorter, whether the image is nearer to the inferior line, superior line, or the center. To define a precise separation between these categories is not possible, but it must be done anyway, and there is where the ambiguity lies.

### 3.3 Data Collection & Methodology for Data Analysis

This questionnaire was applied in the platform 'Moodle institutional' to 170 engineer students, 90 of which study at a public university and the other 80 are students in a private university, both in Chile. The questionnaire was applied between the end of September and the beginning of October 2020, when the students went to online classes due to the pandemic. Every student was connected to the internet and had access to the institutional platform to respond.

From all the students that participated were only analyzed those who completed at least two attempts and that, also, explained how they solved the task. The detail of students that fulfilled these criteria is shown at Table 1.

Table 1: Participants that were analyzed in this research

| Version with Decimal Numbers | Version with Integer Numbers | Total |
|---|---|---|
| 47 | 40 | 87 |

To analyze the data, an analysis of content of the description of the strategies used by the students was made, with the intention of describing the text and extracting inferences of it when relating it with other variables of the research (Mayring, 2015).

The units of analysis are:

- Version of the task: With integer or decimal numbers in the coefficients that define the linear equation which its graph is shown in the task. This variable is defined in the algorithm of the task.
- Question's score: Correct or incorrect. This score is made automatically by the platform.
- Used strategy: From the previously defined theorical elements are deductively raised (Mayring, 2015, p. 374) the following categories: 1) Estimation of the image, which will be labeled as 'to estimate', 2) The calculus of the linear equation and evaluation of the preimage on it, which is labeled as 'to calculate the equation', 3) The use of the two previously listed strategies, e.g., by calculating the linear equation and using the estimation to control the answer, and 4) Some other method that has not been mentioned, which will be labeled as 'other', 5) From the delivered data it is not possible to infer the strategy used by the student, this will be labeled as 'no info'. These categories are disjunct, and they contemplate every possible solution. To identify to which category does each of the students' answers belongs, the shape reduction procedure is used as the technique, to maintain the essential content and to create, by the abstraction, a general vision of the strategies, but getting an image of it (Mayring, 2015, p. 373).

According to the available data, its statistical analysis allows as to assure that the samples (first and second attempt) even when there is a relationship, when using the test of Kolomogorv-Smirnov these does not possess a normal distribution. This is the reason it does not make sense to compare the means. Moreover, when applying the Wilcoxon test to the same set of data we see that the samples are significantly different (at 95% confidence level). That is, the change produced from the first to the second attempt is what in statistics is called 'significant'. Now, by descriptive analysis of the samples, we ensure that 96% of the students that used the strategy 'to estimate' on the first attempt they kept it on the second one. Meanwhile, only 50% of the students who chose 'to calculate the equation' kept the strategy on the second attempt. If we observe what happens with the students that firstly tried to calculate the equation, 46% of them changed their strategy to estimation. The contingency matrix for these two states is:

Table 2: Contingency matrix of the two principal strategies.

| Attempt 2/Attempt 1 | Estimates | Calculates |
|---|---|---|
| Estimates | **96%** | **46%** |
| Calculates | **4%** | **50%** |

After making a content structural analysis from the descriptive statistics, some particular cases are analyzed under these criteria to show some particular elements of the studied phenomenon.

## 4 Results

### 4.1 General Results

When categorizing the strategies according to version and the attempts the results obtained are as synthetized at Table 3 and Figure 6.

According to the data on Table 1 and Figure 6, at the first attempt, the strategies 'to estimate' and 'to calculate the equation' have identical frequency for those students solving the version with integer numbers. However, in the version with decimal numbers, the 'to estimate' strategy appears almost four times more than the 'to calculate the equation' strategy. Based on these data we can conclude that the nature of the coefficients does influence the way the students visualize and solve the task.

Table 3: Used strategy at each attempt according to the task's version

| Strategy | Decimals Numbers First attempt | Decimals Numbers Second Attempt | Integer Numbers First Attempt | Integer Numbers Second Attempt |
|---|---|---|---|---|
| estimates | 66,0% | 74,5% | 45,0% | 75,0% |
| calculates equation | 17,0% | 17,0% | 45,0% | 20,0% |
| both | 2,1% | 2,1% | | |
| visualizes the function | 2,1% | | | |
| Does not remember | 2,1% | | 2,5% | |
| No info. | 10,6% | 6,4% | 7,5% | 5,0% |
| Total | 100% (n=47) | 100% (n=47) | 100% (n=40) | 100% (n=40) |

#### 4.1.1 First attempt's strategy, according to straight line's coefficients' nature

The other strategies appear marginally. In particular, it calls for attention that the strategies 'to estimate' and 'to calculate the equation' appears in a complementary way very few times, that is, to calculate the linear equation and then control the solution from the graph was not an observed strategy.

One of the strategies that showed up unexpectedly was to visualize the function, which consists in visualizing the function as a whole, specifically as the identity function and from that function obtaining the answer.

From the mathematical working space, the strategy 'to calculate the equation' is interpreted as a work that is on the semiotic-instrumental plane, because a visualization is unfolded but it is oriented to the use of the formula of the linear equation as a symbolic artefact. In this work, the instrumental work is so strong that the direct visualization is practically not used to control this algebraical solution, due to the marginal quantity of students that use both strategies jointly.

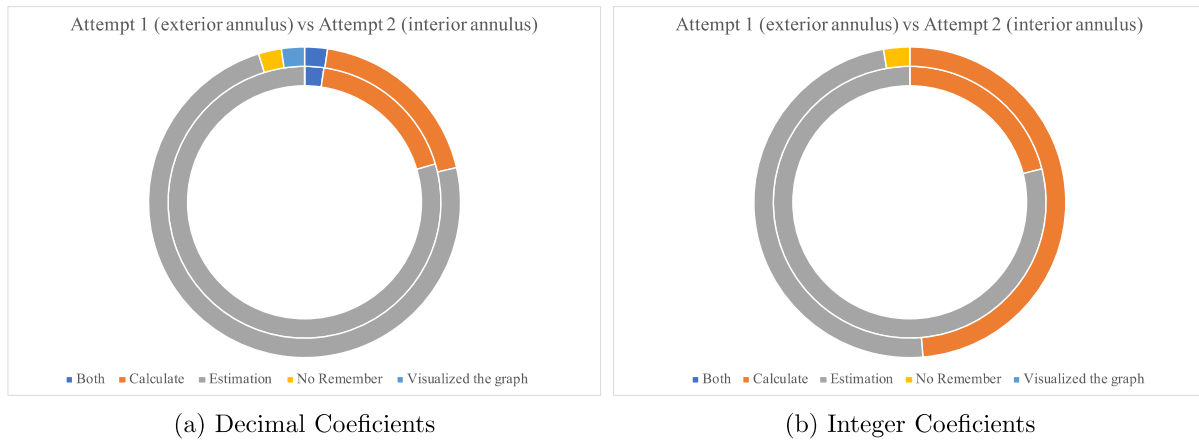

(a) Decimal Coeficients    (b) Integer Coeficients

Figure 6: Graph with the strategies used by the students in attempts 1 and 2 according to the task's version (integer or decimal numbers)

### 4.1.2 Strategy changes for second attempt based on straight line's coefficients' nature influenced by feedback

Once the students answer the first attempt they had access to the feedback, which influenced the strategy used on the second attempt.

On the second attempt of students solving the version with integer numbers, the 'to estimate' strategy's frequency grows in 30%. On the other hand, the 'to calculate the equation' strategy's frequency decreases in a 25%. The marginal strategies disappear, which explains the 5% of difference between the rise of the estimation and the fall of the 'calculating' strategy. In this case, the feedback has a strong influence in the change of mathematical work made by the students, due to their change from a semiotic-instrumental work (when calculating the linear equation) to a purely semiotic work, where the lines of the grid work as semiotic tools that guide the visualization.

In the case of the students working with decimal numbers, the strategy 'to estimate' grows only a 9% its frequency. On the contrary, the strategy 'to calculate the equation' suffers no change, that is to say, the same number of students that used this strategy at the first attempt used it on the second one. The rise of the strategy 'to estimate' is explained by the fall of the other strategies. In this case, the feedback has a more moderated influence on the change of mathematical work, which could be because the possible margin of change is lower, but the striking thing is that those who made a semiotic-instrumental work on the first attempt by working the linear equation, are still doing it on the second attempt. The feedback had no effect on their mathematical work.

### 4.1.3 Strategies effectiveness

Other question that we were able to answer from the data is: What is the effectiveness of the strategies used. From Table 3 and Figure 7 we observe that the strategy 'to estimate' is more effective than 'to calculate the equation'.

This happens on the version with decimal parameters, where we can see that a bit more than half of the students that estimated the image answered correctly, while only a bit more than the 30% of the students that calculated the equation answered correctly, on the same version of the task. This could be because those who calculated the equation had to estimate values that are, by nature, not integers, which increases the error possibilities.

Table 4: Effectiveness of strategy according to version and attempt

| **Decimal Numbers** | correct | incorrect | N (100%) | Difference |
|---|---|---|---|---|
| Estimates | 56,06% | 43,94% | 66 | 12,12% |
| Calculates Equation | 31,25% | 68,75% | 16 | -37,50% |
| **Integer Numbers** | correct | incorrect | N (100%) | Difference |
| Estimates | 70,83% | 29,17% | 48 | 41,67% |
| Calculates Equation | 50,00% | 50,00% | 26 | 0,00% |

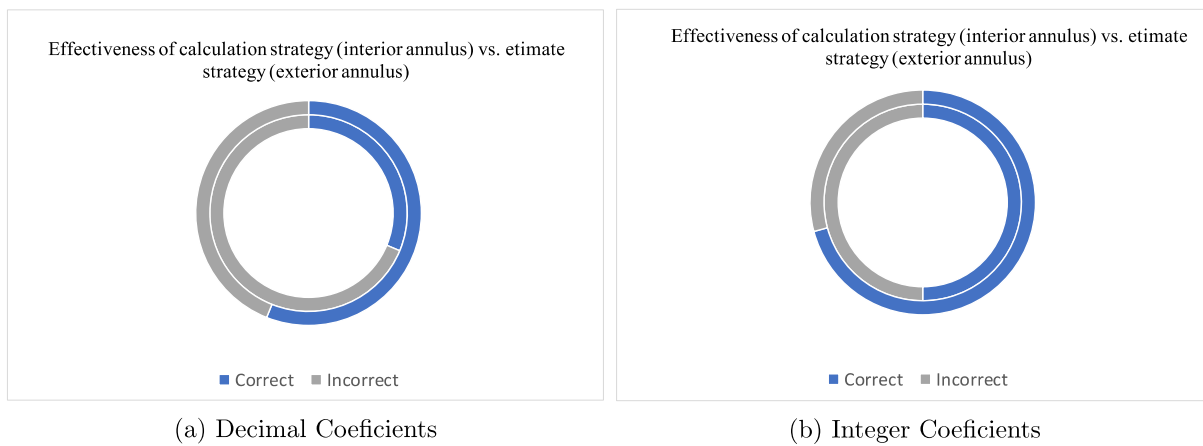

(a) Decimal Coeficients  (b) Integer Coeficients

Figure 7: Graph with the effectiveness of each strategy according to version.

On the other hand, the version containing integer numbers as parameters turned out to be simpler for the students, which is evidenced by the correct answering of more than 70% of students who chose to estimate the image, and by half of the students that chose the strategy 'to calculate the equation' that answered correctly. This version turned out to be less difficult, nevertheless, the estimation is still a more effective strategy, regardless of the version of the task.

### 4.2 Special Results

One of the selected cases belongs to a student that faced the version of the task containing decimal numbers. The task was to estimate the image of x=6 and the student entered y=-23.5 as the answer. The question is shown on Figure 8, along with the answer and the response about the used strategy.

It is seen, in the strategy, that the student makes estimations of the images of x=-3 and x=1 towards integer numbers.

In the case of x=-3, the student approximates the image to y=10. By observing the graph, it is possible to note that the image is close to 10, but it is less than the mentioned value. On the other hand, the student approximates to y=-5 the image of x=1, image that is, indeed, close to

–5, but less than it. Then, when evaluating in the equation, and entering the answer to the system, this considers it as incorrect. What calls the attention of this strategy is that, to calculate the linear equation, the student must make estimations for two points. As seen, in this case the student executes this procedure correctly, although, clearly, there is an inclination toward the integer numbers but, because of the nature of the coefficients, the student commits an error that leads to later errors, enlarging the distance between the obtained image from what the system considers as a good estimation. In some way, the semiotic-instrumental work, oriented by the use of the linear equation is so strong that it changes the way of 'visualizing' the graph, by orienting estimations toward integer values.

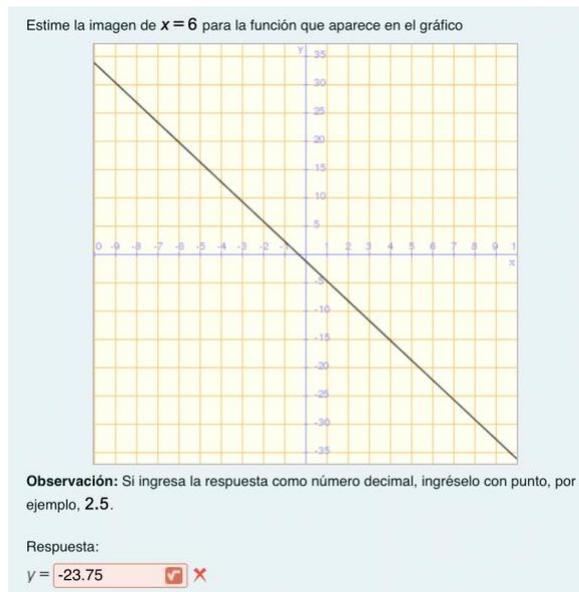

I took as reference two points on the graph: (-3, 10) y (1, -5), to calculate the line's slope m=-15/4 for this way to replace the values in the general equation of the line with the point (1, -5): (y+5) =-15/4(x-1). Finally, I replaced the given value of 'x' and obtained the image of that point on 'y'.

Figure 8: Answer and strategy described by one of the students.

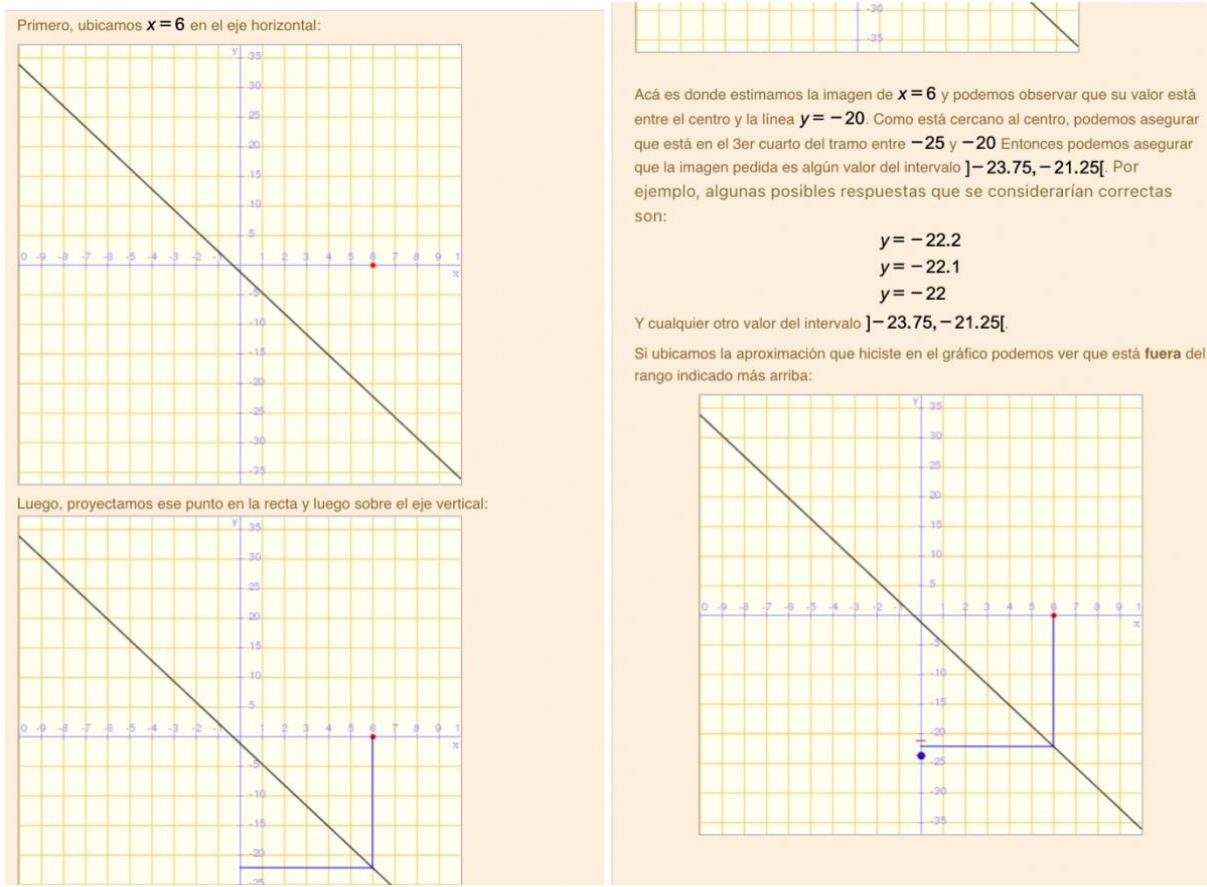
Figure 9: Feedback, delivered by the system, given to the student after answering the task.

When looking at the system's proposed feedback on Figure 9, for one thing it is seen that the answer is on the limit defined by the system as correct. This can be considered as a limitation of the task, particularly of the validation and its associated punctuation. For example, the system could give a partial punctuation or include the extreme values of the interval although, in this case, it would be assumed that the image could be at the start of the first or fourth quarter, more than at the center, which would give a different sense to the visualization of the image. This case shows that there is a level of ambiguity when defining a question like this one, which is proper of the estimation.

The feedback establishes that: "As it is close to the center, we can ensure that it is on the third quarter of the stretch between –25 and –20. Then, we can assure that the image asked is some value of the interval ]-23.75, -21.25[."

Other selected case is that of the student who uses both strategies, meaning that the student coordinates a semiotic-instrumental work, where the linear equation's formula is the symbolic artefact, with a semiotic work that controls the algebraically obtained solution. The interesting thing is that this was a marginal strategy.

"Firstly, I assumed that the slope was –1 (because it did not pass through the origin, but a little below it), this way I tried to calculate the number that multiplies 'x' with the value of x=6, supposing that its image is exactly 2 and, once I had my supposed equation "-4x-1", I replaced the 'x' with 7, getting 29 as the result, which is coherent to what I estimated visually."

Finally, we highlight, from the answers observed on Table 5 in which students describe their work, that estimation is a second-rated strategy, due to students using it only because they did not remember how to do it any other way, moreover, they considered its validity as limited. They tried to calculate the linear equation and, because this was impossible, they estimated. In other words, the purely semiotic work is not considered, by some students, as a mathematical work.

Table 5: Students-made descriptions, which shows that the estimation is a second rated strategy.

| |
|---|
| I don't remember how to do it; I got my result by a visual approximation |
| I tried to find a corner where the line passed exactly to have an idea of the line's slope. But it turns out that was the precise number that I was asked for, according to my suppositions, so I had to do no calculus. |
| Just visually, because the point p (7, -12) is given, and trying to make a function with the intersection with the y-axis would not be exact. |
| First, I tried to find a strategic point, where the line passes through a "pretty" point of the grid. But I could not find one, so I simply used the visual criteria to estimate the image of x=-7 |
| Simply by looking at the graph, because using the point-slope formula gives as a result something quite different to that expressed on the graph. |
| Honestly, I did not know how to do it, I put a number visually. |

## 5  CONCLUSION

In this article, we propose three main questions. The first is: What is the mathematical work that students do when reading a graph if the coefficients that describe an affine function change between whole numbers and decimals? Based on the data collected, we can answer that this work is purely semiotic. By using the grid as a semiotic tool, this can be seen more clearly when the coefficients are decimal numbers. On the other hand, when the coefficients are integers, we see a strong appearance of a semiotic-instrumental work, in which the linear equation is the symbolic artefact used. It is important to note that if the function is changed, for example a quadratic function, it is possible that a smaller portion of the students will use the algebraic strategy. Due to the need for three points and the emergence of a three-variable system of equations, this method would be much harder and could involve knowledge that is not as internalized by students as the linear equation is.

Given the second question: Does a directly value-reading oriented feedback have any influence on the mathematical work made by the students? The answer is yes, but the influence changes according to the task version. The effect is much higher if the coefficients are integers, due to the much lower change margin present, compared to the decimal number task. This calls the attention because, in this version of the task, those who chose calculating the linear equation as a solving strategy did not change said strategy in the second attempt.

The interesting fact is that the final work made by the students is similar despite the task version. In both cases 3 out of 4 students make a direct visualization work, which has certain consequences on the later assessment. If the goal is to measure a direct estimation, without using the linear equation, then it seems more appropriate to work with decimal numbers in the task's coefficients.

Finally, answering the third question: Which is the effectiveness of the work made by the students? We can conclude that, in both cases, the effectiveness of the estimation is higher than the one coming from calculating the linear equation, although (generally), the task containing decimal numbers seems more challenging for the students.

To answer these questions allows us to contemplate the design of the tasks from a wider perspective. Specifically, in this task, it is shown the effect of an implicit variable, as is the type of number defining the line and its respective relation with the grid. If, for example, the challenge is to find the image of a decimal number, it is possible that the mathematical work is also modified, or that the grid is eliminated, which would cause a change in the algorithm to give the correct or incorrect ranges. This allows us to think in new investigations to deeper comprehend the mathematical work based on similar tasks as the one described above, which are very frequent in educational systems.